\documentclass[12pt]{article}
\usepackage{amsmath,amssymb,amsthm,amsfonts}

\theoremstyle{definition}

\theoremstyle{remark}

\numberwithin{equation}{section}

\newtheorem{theorem}{\sc Theorem}[section]
\newtheorem{lemma}[theorem]{\sc Lemma}
\newtheorem{corollary}[theorem]{\sc Corollary}

\newcommand{\R}{\mathbb{R}}
\renewcommand{\Im}{\operatorname{Im}}
\begin{document}

\title {Regularity  of the solutions for nonlinear biharmonic
equations in $\R^N$
\thanks{Research Supported in part by the National Natural Science Foundation of
China (10471052) and the XiaoXiang Funds of Hunan Normal
University.}}
\author{ Yinbin Deng \\
Department of Mathematics,
Huazhong  Normal University\\
 Yi Li\\
  Department of Mathematics, University of Iowa, USA}
\date{}

\maketitle

\begin{abstract}
The purpose of this paper is to establish the regularity the weak
solutions for the nonlinear biharmonic equation
$$
\left\{
\begin{array}{lr}
\Delta^2u+a(x)u=g(x,u),\\
u\in H^2 ({\R^N}),
\end{array}
\right.
$$
where the condition $u\in H^2(\R^N)$ plays the role of a boundary
value condition, and as well expresses explicitly
that the differential equation is to be satisfied in the weak sense.


AMS: 35B40, 35B65, 35G20, 35J60, 35J30

\end{abstract}

\section{\label{S1Int}Introduction}

The purpose of this paper is to establish the regularity  of the
weak solutions for a certain nonlinear biharmonic equation in ${
\R^N}$.  We consider solutions $u\colon {\R^N}\rightarrow {\R}$ of
the problem
$$
\left\{
\begin{array}{lr}
\Delta^2u+a(x)u=g(x,u),\\
u\in H^2 ({\R^N}),
\end{array}
\right.
\eqno{(1.1)}
$$
where the condition $u\in H^2(\R^N)$ plays the role of a boundary
value condition, and as well expresses explicitly
that the differential equation is to be satisfied in the weak sense.
We assume
that
\begin{enumerate}
\raggedright
\item[$H_1)$] $g(x,u)\colon {\R^N}\times {\R^1}\rightarrow {\R^1}
$ is measurable in $x$ and
continuous in $u$, and
$\sup\limits_{\substack{x\in {\R^N}\\0\le u\le M_{\mathstrut}}}
\left|g(x,u)\right|<\infty$
for every $M>0$;

\item[$H_2)$] there exist two constants $\sigma>\delta>0$ and two functions
$b_1(x),b_2(x)\in L^{\infty}(\R^N)$
such that $\left|g(x,u)\right|\le b_1(x)|u|^{\delta + 1}+b_2(x)|u|^{\sigma+1}$;

\item[$H_3)$] $\lim\limits_{\left|x\right|\rightarrow\infty}a(x)=k^2>0$ with $k>0$ and
$(k^2-a(x))\in L^2(\R^N)\cap L^{\infty}(\R^{N})$.
\end{enumerate}
Then we have the following theorems:

\begin{theorem}
\label{Theorem1.1} Assume that $H_1)$ to $H_3)$ hold with
$\sigma+1 <\frac{N+4}{N-4}\, $ if $N\ge 5$. Let $u$ be a weak
solution of \textup{(1.1)}. Then $u\in H^4(\R^N)\cap W^{2,
s}{(\R^N)}$ for $2\le s \le +\infty$. In particular $u\in
C^2(\R^N)$ and
$\lim\limits_{\left|x\right|\rightarrow\infty}u(x)=0$,
$\lim\limits_{\left|x\right|\rightarrow\infty}\Delta u(x)=0$.
\end{theorem}

    Dealing with regularity of solutions
    is much more complicated for biharmonic equations
    than for problems that can be treated by well-developed
    standard methods, such as second-order elliptic problems.
First of all, there is no maximum principle for the biharmonic
problem. So we can't get some estimates of the solutions by the
methods used to deal with second-order elliptic problems.
Secondly, we know little about the properties of the
eigenfunctions of the biharmonic operator in $\R^N$. To overcome
these difficulties, we first introduce the fundamental solutions
for the linear biharmonic operator $\Delta ^2 +k^2 $ for $k >0$.
By applying some properties of Hankel functions, which are the
solutions of Bessel's equation, we obtain the asymptotic
representation of the fundamental solution of $\Delta^2 +k^2 $ at
$\infty$ and $0$. Then we prove that, for $p>1$,
$$
\Delta ^2 -\lambda \colon   \quad  W^{2, p}(\R^N) \longrightarrow
L^p(\R^N)
$$
is an isomorphism if $\lambda <0$. Some estimates of the solutions
of (1.1) can be obtained from the properties of the fundamental
solutions of $\Delta ^2 -\lambda$. We also establish some $L^p$
theory for the biharmonic problem (1.1) so that a bootstrap
argument can be used to deduce the regularity of the solutions of
the biharmonic problem (1.1). Please refer to Grunau \cite{g},
Jannelli \cite{j}, Noussair, Swanson and Yang \cite{nsy}, Peletier
and Van der Vorst \cite{pv}, Pucci and Serrin \cite {ps}
 for the early results
on the existence and other properties of solutions associated with biharmonic operators.

The organization of this paper is as follows:
In Section \ref{S2Fun}, we introduce the fundamental
solutions of $\Delta ^2 -\lambda$ for $\lambda <0$ and establish some properties of
these fundamental solutions.
In Section \ref{S3H4r}, we show that  a weak
solution of the linear problem
$$
\left\{
\begin{array}{lr}
\Delta^2u-\lambda u=f(x),\\
u\in H^2 (\R^N),
\end{array}
\right.
\eqno{(1.2)}
$$
belongs to $H^4(\R^N)$ whenever $f\in L^2(\R^N)$.
In Section \ref{S4W2p}, we obtain a sharper relationship between the
 regularity of the weak solutions of the linear biharmonic problem (1.2) and the
 properties of the inhomogeneous term $f$ in (1.2).
In Section \ref{S5Reg}, we establish the regularity of
the weak solutions for the nonlinear problem (1.1).

\section{\label{S2Fun}Fundamental solutions }

In this section, we give some properties of the fundamental
solutions for the biharmonic operators $\Delta^2+k^2$. The proof
of these properties can be find in \cite {dengli}.

\begin{lemma}
\label{Lemma2.2}\ Let  $G^{(N)}_{k}(|x|)$
 be the fundamental solutions of biharmonic operator $\triangle^2  +k^2$ for $k>0$ and
 $g^{(N)}_{\delta}(|x|)$
 be the fundamental solutions of Laplace operator $-\triangle
 +\delta$. Then we have

\begin{enumerate}
\item[{\textup{i)}}]
$$
G^{(N)}_k(x)\in C^{\infty}({\R^N})\setminus
\{0\}
$$
and
$$
\Delta^2 G^{(N)}_k(x)+k^2G^{(N)}_k(x)=0  {\text{\qquad for }}  x\not=0\,.
\eqno{(2.20)}
$$

\item[{\textup{ii)}}]  As $\left|x\right|\rightarrow\infty$,
$$
e^{(\sqrt{k}/\sqrt{2})\left|x\right|}G^{(N)}_k(x)\rightarrow 0
{\text{\quad and\quad }}  e^{(\sqrt{k}/\sqrt{2})\left|x\right|}
\left|\nabla G^{(N)}_k(x)\right|
\rightarrow 0\,.
\eqno{(2.21)}
$$

\item[{\textup{iii)}}]   As $\left|x\right|\rightarrow 0$,
\begin{align*}
&
\begin{aligned}
G^{(N)}_k(r)
&=\frac{2^{\nu-2}\Gamma(\nu-1)}{2(2\pi)^{N/2}}r^{2-2\nu}\\
&\qquad+
O(r^{4-2\nu})
\vphantom{\frac{2^{\nu-2}\Gamma(\nu-1)}{2(2\pi)^{N/2}}}
\end{aligned}
& &
\text{\quad if }\nu=\frac{N-2}{2}>1\text{ and }\nu\notin {\mathcal{N}};
\\[3\jot]
&
\begin{aligned}
G^{(N)}_k(r)
&=\frac{2^{\nu-2}\Gamma(\nu-1)}{2(2\pi)^{N/2}}r^{2-2\nu}\\
&\qquad+
O(r^{4-2\nu}+{\ln r})
\vphantom{\frac{2^{\nu-2}\Gamma(\nu-1)}{2(2\pi)^{N/2}}}
\end{aligned}
& &
\text{\quad if }\nu=\frac{N-2}{2}\ge 2\text{ and }\nu\in {\mathcal{N}};
\\[3\jot]
&
\begin{aligned}
G^{(N)}_k (r)
&\approx O(\ln r)
\end{aligned}
& &
\text{\quad if } N=4 \ (\textstyle\nu=\frac{N-2}{2}=1 );
\\[3\jot]
&
\begin{aligned}
G^{(N)}_k (r)
&=O(1)
\end{aligned}
& &
\text{\quad if } N=2,3 \ (\textstyle\nu=0,\frac{1}{2}).
\end{align*}

\item[{\textup{iv)}}]   $|G^{(N)}_k(r)|\le Cg^{(N)}_{\delta}(r)$
for some positive constants $C$ and
$0<\delta<\frac{\sqrt{k}}{\sqrt{2}}$.
\end{enumerate}
\end{lemma}

It
follows from properties {\textup{ii)}} and {\textup{iii)}} of Lemma \ref{Lemma2.2}
that:

\begin{corollary}
\label{Corollary2.3}
$$
\begin{aligned}
&G^{(N)}_k(x)\in L^p({\R^N}) & & \text{ for }  \ 1\le p\le +\infty , & & \text{ if } \ N=2,3,\\
&G^{(N)}_k(x)\in L^p({\R^N}) & & \text{ for }  \ 1\le p<+\infty , & & \text{ if } \ N=4,\\
&G^{(N)}_k(x)\in L^p({\R^N}) & & \text{ for }  \ 1\le p<\tfrac{N}{N-4}, & & \text{ if } \ N\ge 5,\\
&\left|\nabla G^{(N)}_k(x)\right|\in L^p  & & \text{ for } \ 1\le p<\tfrac{N}{N-3}, & & \text{ if } \ N>3,\\
&\left|\nabla G^{(N)}_k(x)\right|\in L^p  & & \text{ for } \ 1\le p<+\infty, & & \text{ if } \ N=3,\\
&\left|\nabla G^{(N)}_k(x)\right|\in L^p  & & \text{ for } \ 1\le p\le +\infty, & & \text{ if } \ N=2,\\
&\left|\Delta G^{(N)}_k(x)\right|\in L^p  & & \text{ for } \ 1\le p<\tfrac{N}{N-2}, & & \text{ if } \ N\ge 3,\\
&\left|\Delta G^{(N)}_k(x)\right|\in L^p  & & \text{ for } \ 1\le p<+\infty, & & \text{ if } \ N=2.
\end {aligned}
\eqno {(2.24)}
$$
\end{corollary}

Using this information about $G^{(N)}_k(x)$, we can express
solutions of the inhomogeneous biharmonic equation as convolutions
of fundamental solutions with the inhomogeneous term. The
following Theorem can also be found in \cite {dengli}.

\begin{theorem}
\label{Theorem2.4}\

\begin{enumerate}
\item[{\textup{i)}}]  Let $f\in L^2({\R^N})\cap L^{\infty}({\R^N})$
and
$$
u=\int_{\R^N}f(z)G^{(N)}_k(x-z)\,dz\,.
$$
Then
$$
\Delta^2u+k^2u=f(x)\,.
$$

\item[{\textup{ii)}}] Let $u$ be a distribution such that
$$
\Delta^2 u+k^2u=f
$$
and $f\in L^2({\R^N})\cap L^{\infty}({\R^N})\,$. Then
$$
u=\int_{\R^N}f(z)G^{(N)}_k(x-z)\,dz\,.
\eqno{(2.25)}
$$

\item[{\textup{iii)}}] There are no nontrival distributions such that
$$
\left\{
\begin{array}{cl}
\Delta^2u+k^2u=0\,,\\
u\in W^{2,2}({\R^{N}})\,.
\end{array}
\right. \eqno{(2.26)}
$$
\end{enumerate}
\end{theorem}

\section{\label{S3H4r}$H^4$-regularity}

The main purpose of this section is to show that  a weak
solution of the linear problem
$$
\left\{
\begin{array}{lr}
\Delta^2u-\lambda u=f(x),\\
u\in H^2 (\R^N),
\end{array}
\right.
\eqno{(3.1)}
$$
belongs to $H^4(\R^N)$ whenever $f\in L^2(\R^N)$. To this end, we recall a
well-known result which can be found in \cite {rs1}.

\begin{lemma}
\label{Lemma3.1} Let $h\in L^p(\R^N)$ for some $p\in [1,+\infty]$,
and consider the equation
$$-\Delta u+u=h\eqno{(3.2)}$$
in the sense of distributions.

\begin{enumerate}
\item[{\textup{a)}}]  There is a unique tempered distribution $u=\Gamma(h)$
satisfying \textup{(3.2)}.

\item[{\textup{b)}}]  If $h\in L^p(\R^N)$ for some $p\in (1,+\infty)$, then
$\Gamma (h)\in W^{2,p}(\R^N)$ and there exists a constant $C(N,p)$ such that
$$
\left\|\Gamma (h)\right\|_{W^{2,p}}\le C(N,p)\left\|h\right\|_{L^p}
$$
for all $h\in L^p(\R^N)$.

\item[{\textup{c)}}] For $p\in (1,+\infty)$,
$-\Delta+1\colon W^{2,p}({\R^N})\rightarrow
L^p(\R^N)$ is an isomorphism.
\end{enumerate}
\end{lemma}

By applying this lemma, we can obtain the $W^{4,p}(\R^N)$ regularity
for the linear biharmonic problem.

\begin{lemma}
\label{Lemma3.2} Let $v\in W^{2,p}(\R^N)$, $ w\in L^p(\R^N)$
for some $p\in (1,+\infty)$ be such that
$$\int_{\R^N}\Delta v\Delta z\,dx = \int_{\R^N}wz\,dx
\text{\qquad for all} \ z\in C^{\infty}_0(\R^N)\,.\eqno{(3.3)}$$
Then $v\in W^{4,p}(\R^N)$ and $\Delta^2 v=w$.
\end{lemma}

\begin{proof}
{}From (3.3), it follows that $u=\Delta v$ is a distribution solution of
$$\Delta u=w \ \ {\text{ and }} \ \
u,w\in L^p(\R^N)\,.\eqno{(3.4)}$$
Thus
$$(-\Delta+1)u=u-w\in L^p(\R^N)\,.$$
By applying Lemma \ref{Lemma3.1} we find
that
$-\Delta+1 \colon W^{2,p}({\R^N})\rightarrow L^p(\R^N)$ is an isomorphism. So there
 exists $\varphi\in W^{2,p}(\R^N)$ such that
$$(-\Delta+1)\varphi=u-w\,,$$
that is
$$-\int_{\R^N}\varphi\Delta z\,dx+\int_{\R^N}\varphi z\,dx=\int_{\R^N}uzdx-\int_{\R^N}wz\,dx$$
for all $z\in C^{\infty}_0{(\R^N)}$. From (3.3) we have
$$\int_{\R^N}(\varphi-u)\Delta z\,dx=\int_{\R^N}(\varphi-u)z\,dx
\text{\qquad for all} \ z\in C^{\infty}_0(\R^N)$$ and hence
$$\int_{\R^N}(\varphi-u)(-\Delta z+z)\,dx=0
\text{\qquad for all} \ z\in C^{\infty}_0(\R^N)\,.\eqno{(3.5)}$$
Consider the equation
$$-\Delta z+z=|\varphi-u|^{p-2}(\varphi-u)\,.\eqno{(3.6)}$$
It follows from $\varphi-u\in L^p$ that
$|\varphi-u|^{p-2}(\varphi-u)\in L^{p'}(\R^N)$ with
$\frac{1}{p}+\frac{1}{p'}=1$. By Lemma \ref{Lemma3.1},
the problem (3.6)
possesses a unique solution $z\in W^{2,p}(\R^N)$. Since
$C^{\infty}_0(\R^N)$ is dense in $W^{2,p'} (\R^N)$, we can find a
sequence $ \{{z_n} \} \subset C^{\infty}_0 (\R^N)$ such that
$$z_n\rightarrow z  {\text{\quad in }}  W^{2,p'}(\R^N)  {\text{\quad as }}
  n \rightarrow\infty.$$
{}From (3.5) and (3.6), we have
\begin{align*}
0=\int_{\R^N}(\varphi-u)(-\Delta z_n+z_n)\,dx
&\rightarrow\int_{\R^N}(\varphi-u)(-\Delta z+z)\,dx\\
&=\int_{\R^N}|\varphi-u|^p \,dx\,.
\end{align*}
This implies that $\varphi-u\equiv 0$ and hence $u\in W^{2,p}(\R^N)$. It follows
that $v\in W^{4,p}(\R^N)$.
\end{proof}

To obtain the $H^4$-regularity of solutions of (3.1), we rewrite
the
problem (3.1) in the form
$$
\left\{
\begin{array}{lr}
\Delta^2u=f+\lambda u,\\
u\in H^2(\R^N).
\end{array}
\right.
$$
The $H^4$-regularity of solutions of (3.1) follows from Lemma \ref{Lemma3.2}.
In fact, we can
get
a
more general result:

\begin{lemma}
\label{Lemma3.3} Let $f\in L^p(\R^N)$ for some $p\in (1,+\infty)$,
and let
 $u$ be the solution of
the
 linear biharmonic problem
$$
\left\{
\begin{array}{lr}
(\Delta^2-\lambda)u=f,\\
u\in W^{2,p}(\R^N).
\end{array}
\right.
\eqno{(3.7)}
$$
Then $u\in W^{4,p}(\R^N)$.
\end{lemma}

In the following lemma, we show that
the
problem (3.7) possesses a unique solution
$u\in W^{2,p} (\R^N)$ for given
$p\in [2,+\infty)$ if $f\in L^p(\R^N)$ and $\lambda<0$.

\begin{lemma}
\label{Lemma3.4} For $p\in [2,+\infty)$, $f\in L^p(\R^N)$,
the
problem \textup{(3.7)} possesses a unique solution if $\lambda<0$.
\end{lemma}

\begin{proof}
{}From Lemma \ref{Lemma3.3}, the solution of (3.7) must belong to
$W^{4,p}(\R^N)$.
 Suppose that $u\in W^{4,p}({\R^N})$
 is a solution of the homogeneous problem
$$
\left\{
\begin{array}{lr}
\Delta^2u-\lambda u=0,\\
u\in W^{4,p}({\R^N}).
\end{array}
\right. \eqno{(3.8)}
$$
Rewrite (3.8) in the form
$$
\left\{
\begin{array}{lr}
-\Delta(-\Delta u)=\lambda u,\\
u\in W^{4,p}({\R^N}),\;-\Delta u\in W^{2,p}({\R^N}).
\end{array}
\right.
$$
By using a bootstrap argument, it follows that
$$u\in C^4 ({\R^N})\cap L^p({\R^N}),
 \qquad \Delta u\in C^2({\R^N})\cap L^p({\R^N}),$$
and $\lim\limits_{\left|x\right|\rightarrow\infty}u(x)=0$,  $
\lim\limits_{\left|x\right|\rightarrow\infty}\Delta u(x)=0$. Define
$u_1=(\Delta-\sqrt{\lambda})u$, $u_2=(\Delta+\sqrt{\lambda})u$.
Then
$$(\Delta+\sqrt{\lambda})u_1=0, \qquad (\Delta-\sqrt{\lambda})u_2=0,\eqno{(3.9)}$$
and $$u=\frac{1}{2\sqrt{\lambda}}(u_2-u_1),$$
$$\lim_{\left|x\right|\rightarrow\infty}u_1(x)=0 , \qquad
\lim_{\left|x\right|\rightarrow\infty}u_2(x)=0.
$$
For
$\lambda<0$, the solution of (3.9) can be expressed
in terms of
Hankel
functions.  By the asymptotic behavior of Hankel functions (see
(2.8)), we can deduce that
$$
e^{\Im (\lambda)^{1/4}\left|x\right|}u_i(x)\rightarrow 0,
\quad i=1,2, \ \
\text{\quad as }  \left|x\right|\rightarrow \infty.$$ Thus we have
$$e^{\Im (\lambda)^{1/4} \left|x\right|}u(x)\rightarrow 0
\text{\qquad as }  \left|x\right|\rightarrow \infty;\eqno{(3.10)}$$
it follows from (3.10) that $u\in L^r({\R^N})$ for all $r\in
[2,+\infty)$.
 In particular, $u\in L^2({\R^N})$ and hence $u\in H^{4}({\R^N})$.
 Theorem  \ref{Theorem2.4} gives us that $u\equiv 0$.

This completes the proof of our lemma.
\end{proof}

\section{\label{S4W2p}$W^{2,p}(\R^N)$-regularity}

In this section, we obtain a sharper relationship between the
 regularity of the weak solutions of
the
 linear biharmonic problem (3.1) and the
 properties of the inhomogeneous term $f$ in (3.1).

Recalling
 the properties of the fundamental solution $G^{(N)}_k$ for $k>0$
(see Corollary \ref{Corollary2.3}), Young's  inequality for
convolutions \cite {T} shows
that the convolution $f*G^{(N)}_k$
defines an element of $L^s(  \R^N)$ subject to the following restrictions:
$$
\left\{
\begin{array}{ll}
p\le s\le +\infty & {\text{ if }} \ p>\frac{N}{4}\,,\\
p\le s < +\infty & {\text{ if }} \ p=\frac{N}{4}\,,\\
p\le s\le\frac{Np}{N-4p} & {\text{ if }} \ 1\le p<\frac{N}{4}\,.
\end{array}
\right.
\eqno{(4.1)}
$$
Setting $T_kf=f*G^{(N)}_k$, we see that
$$
T_k\colon  L^p(\R^N)\rightarrow L^s(\R^N) \text{ is a bounded linear operator.}
\eqno {(4.2)}
$$
Referring again to Corollary \ref{Corollary2.3},
we can deduce that for $i=1,2,\dots, N$,
 the convolution $f*\partial _i G^{(N)}_k$ defines
an element of  $L^s$ whenever $f\in L^p$ subject to the restrictions
$$
\left\{
\begin{array}{ll}
p\le s \le +\infty & {\text{ if }} \ p>\frac{N}{3}\,,\\
p\le s<+\infty & {\text{ if }} \ p=\frac{N}{3}\,,\\
p\le s<\frac{Np}{N-3p} & {\text{ if }} \ 1\le p<\frac{N}{3}\,;
\end{array}
\right.
\eqno{(4.3)}
$$
and the convolution $f*\Delta G^{(N)}_k$ defines an element of $L^s$ whenever $f\in L^p$
subject to the restrictions
$$
\left\{
\begin{array}{ll}
p\le s \le +\infty & {\text{ if }} \ p>\frac{N}{2}\,,\\
p\le s<+\infty & {\text{ if }} \ p=\frac{N}{2}\,,\\
p\le s<\frac{Np}{N-2p} & {\text{ if }} \ 1\le p<\frac{N}{2}\,.
\end{array}
\right.
\eqno{(4.4)}
$$
Setting $S^i_k f=f*\partial_i G^{(N)}_K$, $i=1,2,\dots,N$, and
$S^{\Delta}_k f=f*\Delta G^{(N)}_k$,
we see that
\begin{align*}
S^i_k \colon  L^p\rightarrow L^s
&  {\text{ is  a  bounded  linear  operator}}\tag*{(4.5)}\\
\intertext{under the restrictions (4.3) and}
S^{\Delta}_k \colon  L^p\rightarrow L^s
&  {\text{ is  a  bounded  linear  operator}}\tag*{(4.6)}
\end{align*}
under the restrictions (4.4).

\begin{theorem}
\label{Theorem4.1} Given $k>0$ and $f\in C^2_0(\R^N)$, set
$$T_kf(x)=f*G^{(N)}_k(x)  {\text{\qquad for }}  x\in{\R^N}\,.\eqno{(4.7)}$$
Then $T_k f\in C^4({\R^N})$, $\lim\limits_{\left|x\right|\rightarrow\infty}T_k f(x)=0$, and
 $u=T_k f$ satisfies the biharmonic equation
$$\Delta^2 u=\lambda u+f  {\text{\qquad on }}  {\R^N},\eqno{(4.8)}$$
where $\lambda=-k^2$. Furthermore, for all $x\in {\R^{N}}$, the
following formulae are valid for $i,j,l,m=1,2,\dots,N$:
\begin{align*}
T_kf(x) &=\int\! f(x-z)G^{(N)}_k(z)\,dz
=\int\! G^{(N)}_k(x-z)f(z)\,dz,
\\
\partial_i T_k f(x) &= \int\! \partial_i f(x-z)G^{(N)}_k(z)\,dz
=\int\! \partial_i G^{(N)}_k(x-z)f(z)\,dz,
\\
\partial _j\partial_i T_k f(x)
&=\int\!\partial_i f(x-z)\partial_j G^{(N)}_k(z)\,dz,
\\
\partial _m\partial _j\partial_i T_k f(x)
&=\int\!\partial_m\partial_i f(x-z)\partial_j G^{(N)}_k(z)\,dz
=\int\!\partial_i f(x-z)\partial_m\partial_j G^{(N)}_k (z)\,dz\rlap{,}\!
\\
\partial_l\partial_m\partial_j\partial_i T_k f(x)
&=\int\!\partial_m\partial_j f(x-z)\partial_l\partial_i G^{(N)}_k(z)\,dz.
\end{align*}
\end{theorem}

\begin{proof}
Noting that
$$T_kf(x)=T_1f_k(kx)
 \text{\qquad where }  f_k(y)=k^{-4}f\left(\frac{z}{k}\right),$$
we see that, by a change of scale, it is enough to treat the case $k=1$. In the
following, we take $k=1$ and simplify the notation by setting $T_1=T$,
$G^{(N)}_1=G^{(N)}$. Since $f$, $\partial_i f$ and $\partial_{ij}f\in C^0(\R^N)$
and $\partial_i G^{(N)},\partial_{ij}G^{N}\in L^1(\R^N)$, it
follows that the convolutions
\begin{gather*}
Tf=f*G^{(N)}, \qquad \partial_i f* G^{(N)}, \qquad f*\partial G^{(N)},\\
\partial_i f*\partial_j G^{(N)}, \qquad \partial_{ij}f*\partial_{m}G^{(N)}, \qquad
\partial_i f*\partial_{jm} G^{(N)},
\text{\quad and\quad }\partial_{ij}f*\partial_{ml}G^{(N)}
\end{gather*}
are defined and are continuous on ${\R^N}$. They all tend to zero as
$\left|x\right|\rightarrow\infty$.
Hence to prove the theorem it is sufficient to establish the following statements.
\begin{enumerate}
\item[1)]  $\partial_i Tf$ exists and $\partial_i Tf=\partial_if*G^{(N)}$.

\item[2)]  $\partial_if*G^{(N)}=f* \partial_i G^{(N)}$.

\item[3)]  $\partial_{ij}Tf$ exists and $\partial_j\partial_i Tf
=\partial_jf*\partial_iG^{(N)}$.

\item[4)]  $\partial_{mji}Tf$ exists and $\partial_{mij}Tf
=\partial_m\partial_jf*\partial_iG^{(N)}
=\partial_jf*\partial_m\partial_iG^{(N)}$.

\item[5)]  $\partial_{lmji}Tf$ exists and  $\partial_{lmji}Tf
=\partial_m\partial_jf*\partial_l\partial_iG^{(N)}$.

\item[6)]  $\Delta^2Tf+Tf=f$  on ${\R^N}$.
\end{enumerate}

(1) Let  $e_i$ be an element of the usual basis for ${\R^N}$ and $h$ a
non-zero real number. Then
$$\frac{Tf(x+he_i)-Tf(x)}{h}=\int_{\R^N}\frac{f(x+he_i-z)
-f(x-z)}{h}G^{(N)}(z)\,dz$$
and
$$\lim_{h\rightarrow 0}\frac{f(x+he_i -z)-f(x-z)}{h}=\partial_if(x-z).$$
Also
\begin{align*}
\left|\frac{f(x+he_i-z)-f(x-z)}{h}\right|
&\le\left|\frac{1}{h}\int^1_0\frac{d}{dt}f(x+the_i-z)\,dt\right|\\
&=\left|\int^1_0\partial_i f(x+the_i-z)\,dt\right|\\
&\le \max_{z\in{\R^N}}\left|\partial_i f(z)\right|=
\left|\partial_i f\right|_{\infty}.
\end{align*}
Hence, by the dominated convergence theorem,
$$\lim_{h\rightarrow 0}\frac{Tf(x+he_i)-Tf(x)}{h}
=\int_{\R^N}\partial_i f(x-z)G^{(N)}(z)\,dz.$$

(2) For $i=1,2,\dots,N$,
\begin{align*}
\partial_i f*G^{(N)}(x)
&=\lim_{\epsilon\rightarrow 0}\int_{|z|\ge\epsilon}\partial_if(x-z)G^{(N)}(z)\,dz\\
&=-\lim_{\epsilon\rightarrow 0}\int_{|z|\ge\epsilon}\frac{\partial}{\partial z_i}f(x-z)G^{(N)}(z)\,dz\\
&=\lim_{\epsilon\rightarrow 0}\left\{\int_{|z|=\epsilon}\frac{z_i}{|z|}f(x-z)G^{(N)}(z)\,dz
+\int_{|z|\ge\epsilon}f(x-z)\partial_i G^{(N)}(z)\,dz\right\}.
\end{align*}
Now from Lemma \ref{Lemma2.2},
\begin{align*}
\left|\int_{|z|=\epsilon}\frac{z_i}{|z|}f(x-z)G^{(N)}(z)\,dz\right|
&\le \left|f\right|_{\infty}\int_{|z|=\epsilon} \left|G^{(N)}(z)\right|\,dz\\
&=\left|f\right|_{\infty}\left|G^{(N)}(\epsilon)\right|\int_{|z|=\epsilon} \! z\,dz
=\left|f\right|_{\infty}\left|G^{(N)}(\epsilon)\right|\cdot w_N\epsilon^{N-1}\\
&=\left\{
\begin{aligned}
&
\frac{2^{\nu-2}\Gamma(\nu-1)}{2(2\pi)^{N/2}}\epsilon^{4-N}
 w_N\epsilon^{N-1}\left|f\right|_{\infty}
& & {\text{ if }}  N\ge 5,\\
& \vphantom{\frac{2^{\nu-2}\Gamma(\nu-1)}{2(2\pi)^{N/2}}}
O\left(\left|\ln  \epsilon\right|\right)\epsilon^{N-1}  w_N \left|f\right|_{\infty}
& & {\text{ if }}  N=4,\\
& \vphantom{\frac{2^{\nu-2}\Gamma(\nu-1)}{2(2\pi)^{N/2}}}
O(1)\epsilon^{N-1}  w_N \left|f\right|_{\infty}
& & {\text{ if }}  N=2,3,
\end{aligned}
\right.
\end{align*}
where $\nu=\frac{N-2}{2}$. Hence
$$\lim_{\epsilon\rightarrow 0}\int_{|z|=\epsilon}\frac{z_i}{|z|}f(x-z)G^{(N)}(z)\,dz=0$$
and
$$\partial_if*G^{(N)}(x)=\lim_{\epsilon\rightarrow 0}
\int_{|z|\ge\epsilon}f(x-z)\partial_iG^{(N)}(z)\,dz
=f*\partial_iG^{(N)}.$$

(3)  Repeat the proof of (1) with $G^{(N)}$ replaced by $\partial_i G^{(N)}$.

(4) Repeat the proof of (1) and (2) with $G^{(N)}$ and $f$ replaced by $\partial_iG^{(N)}$ and $\partial_j f$.

(5) Repeat the proof of (1) with $G^{(N)}$ and $f$ replaced by $\partial_{il}G^{(N)}$ and $\partial_j f$.

(6)
\begin{align*}
\Delta^2 Tf(x)&=\int\Delta f(x-z)\Delta G^{(N)}(z)\,dz\\
&=\lim_{\epsilon\rightarrow 0}\int_{|z|\ge \epsilon}
\Delta _z f(x-z)\Delta G^{(N)}(z)\,dz\\
&=\lim_{\epsilon\rightarrow 0}\left\{\int_{|z|=\epsilon}
f(x-z)\frac{\partial\Delta G^{(N)}}{\partial r}-\Delta G^{(N)}
\frac{\partial f(x-z)}{\partial r}\,dz
\vphantom{+\int_{|z|\ge\epsilon}f(x-z)\Delta ^2G^{(N)}(z)\,dz}\right.\\
&\qquad+
\left.\vphantom{\int_{|z|=\epsilon}
f(x-z)\frac{\partial\Delta G^{(N)}}{\partial r}-\Delta G^{(N)}
\frac{\partial f(x-z)}{\partial r}\,dz}
\int_{|z|\ge\epsilon}f(x-z)\Delta ^2G^{(N)}(z)\,dz\right\}\\
&=\lim_{\epsilon\rightarrow 0}\int_{|z|=\epsilon}
\left(f(x-z)\frac{\partial\Delta G^{(N)}}{\partial r}
-\Delta G^{(N)}\frac{\partial f(x-z)}{\partial r}\right)\,dz-Tf(x).
\end{align*}
Since
$$
\Delta ^2G^{(N)}(z)=-G^{(N)}(z)
\text{ for all } z\not=0,
$$
we obtain that
\begin{align*}
\lim_{\epsilon\rightarrow 0}\int_{|z|=\epsilon}
\Delta G^{(N)}\frac{\partial f(x-z)}{\partial r} \,dz
&= 0,\tag*{(4.9)}\\
\lim_{\epsilon\rightarrow 0}\int_{|z|=\epsilon}f(x-z)
\frac{\partial}{\partial r}
\left(\Delta G^{(N)}(z)\right)\,dz
&= f(x).\tag*{(4.10)}
\end{align*}

In fact, from  (2.13),
$$(G^{(N)}(r))'=-2\pi rG^{(N+2)} (r).$$
Thus
\begin{align*}
\left(G^{(N)}(r)\right)''
&= 2^2\pi^2r^2G^{(N+4)}(r)-2\pi G^{(N+2)}(r),\\[3\jot]
\Delta G^{(N)}(x)
&= \left(G^{(N)}(r)\right)''+\frac{N-1}{r}\left(G^{(N)}(r)\right)'\\
&= 4\pi^2r^2G^{(N+4)}(r)-2\pi N G^{(N+2)} (r),\\[3\jot]
\left(\Delta G^{(N)} (r)\right)'_r
&= -16\pi^3r^3G^{(N+6)}(r)+
(8+4N)\pi^2rG^{(N+4)}(r).
\end{align*}
By the asymptotic behavior of
$G^{(N)}_k(r)$ (see Lemma \ref{Lemma2.2}) we deduce that,
 as $r=\left|x\right|\rightarrow 0$,
\begin{align*}
\Delta G^{(N)}(x)
&\approx 4{\pi}^2 r^2 G^{(N+4)}(r)
\approx\frac{2^{(N/2)-1}\Gamma\left(\frac{N}{2}\right)}{2(2\pi)^{N/2}}r^{2-N},
\tag*{(4.11)}\\[3\jot]
\left(\Delta G^{(N)}(r)\right)'_r
&\approx -16 {\pi}^3 r^3 G^{(N+6)}(r)
\approx\frac{2^{N/2}\Gamma\left(\frac{N}{2}+1\right)}{2(2\pi)^{N/2}}r^{1-N}.
\tag*{(4.12)}
\end{align*}
Thus
\begin{align*}
\left|\int_{|z|=\epsilon}\Delta G^{(N)}(z)\frac{\partial f(x-z)}{\partial r}\,dz
\vphantom{\frac{2^{(N/2)-1}\Gamma\left(\frac{N}{2}\right)}{2(2\pi)^{N/2}}}
\right|
&\approx \left|\int_{|z|=\epsilon}
\frac{2^{(N/2)-1}\Gamma\left(\frac{N}{2}\right)}{2(2\pi)^{N/2}}
|z|^{2-N}
\frac{\partial f(x-z)}{\partial r}\,dz\right|\\
& \le \left|\nabla f\right|_{L^{\infty}}
\int_{|z|=\epsilon}\frac{2^{(N/2)-1}\Gamma\left(\frac{N}{2}\right)}{2(2\pi)^{N/2}}
|z|^{2-N}\,dz\\
&= \left|\nabla f\right|_{L^{\infty}}\frac{2^{(N/2)-1}\Gamma\left(\frac{N}{2}\right)}{2(2\pi)^{N/2}}
\epsilon^{2-N}\epsilon^{N-1}w_N\rightarrow 0\\
& \qquad\text{\qquad as }  \epsilon\rightarrow 0.
\end{align*}
This gives (4.9).

Now we are going to prove (4.10). From (4.12) we have
$$\epsilon^{N-1}
\left.\left(\Delta G^{N}(r)\right)'_r\right|_{r=\epsilon}\rightarrow
\frac{2^{N/2}\Gamma\left(\frac{N}{2}+1\right)}{(2\pi)^{N/2}}
\text{\qquad as }  \epsilon\rightarrow 0.$$
Thus
\begin{align*}
\int_{|z|=\epsilon}  &f(x-z)\left(\Delta G^{(N)}(z)\right)'_r\,dz
= \int_{|z|=\epsilon}(f(x-z)-f(x))\left(\Delta G^{(N)}(z)\right)'_r \,dz\\
& \qquad+
f(x)\left.\left(\Delta G^{(N)}(r)\right)'_r\right|_{r=\epsilon}\cdot \epsilon^{N-1} w_N\\
&\approx \int_{|z|=\epsilon}(f(x-z)-f(x))\cdot
\frac{2^{N/2}\Gamma\left(\frac{N}{2}+1\right)}{(2\pi)^{N/2}}r^{1-N}\,dz\\
& \qquad+
f(x)\frac{2^{N/2}\Gamma\left(\frac{N}{2}+1\right)}{(2\pi)^{N/2}}w_N\\
&= \frac{2^{N/2}\Gamma\left(\frac{N}{2}+1\right)}{(2\pi)^{N/2}}\epsilon^{1-N}
\int_{|z|=\epsilon}(f(x-z)-f(x))\,dz\\
& \qquad+
f(x)\cdot\frac{\frac{N}{2}\Gamma\left(\frac{N}{2}\right)}{{\pi}^{N/2}}
\cdot\frac{2\pi^{N/2}}{N\Gamma\left(\frac{N}{2}\right)}\\
&= \frac{2^{N/2}\Gamma\left(\frac{N}{2}+1\right)}{(2\pi)^{N/2}}\epsilon^{1-N}
\int _{|z|=\epsilon}(f(x-z)-f(x))\,dz+f(x).
\end{align*}
The limit
(4.10) follows from the facts that
$$\left|f(x-z)-f(x)\right|\le\left|\nabla f\right|_{L^{\infty}}|z|$$
for all $z\in {\R^N}$
and hence
$$\left|\epsilon^{1-N}\int_{|z|=\epsilon}(f(x-z)-f(x))\,dz\right|
\le \epsilon^{1-N} \left|\nabla f\right|_{L^{\infty}}\epsilon
\cdot\epsilon^{N-1}w_n\rightarrow 0.$$ This completes the proof of
Theorem \ref{Theorem4.1}.
\end{proof}

\begin{theorem}
\label{Theorem4.2} Let $\lambda=-k^2$ where $k>0$ and let $f\in L^p$
where $p\in\left(\frac{2N}{N+4},2\right]$. Then $T_kf\in H^2(\R^N)$ and
$\partial_i T_k f=S^i_k f$ and $\Delta T_k f=S^{\Delta}_k f$.
Furthermore, $T_kf$ is
a weak solution of \textup{(3.1)}, where $S^j_k$ and
$S^{\Delta}_k$ are given by \textup{(4.5)} and \textup{(4.6)}.
\end{theorem}

\begin{proof}
Let $\{f_n\}\subset C^2_0$ be a sequence such that
$|f_n-f|_p\rightarrow 0$ as $n\rightarrow\infty$. Since
$p>\frac{2N}{N+4}$, we have $\frac{2N}{N-4}<\frac{Np}{N-2p}$ when $N>4$ and
$p<\frac{N}{4}$. From (4.2) it follows that
$$
T_kf_n\text{ and }T_kf\in L^s
$$
and that
$$
\left|T_kf_n-T_kf\right|_{L^s}\rightarrow 0\text{ as }
n\rightarrow\infty,
$$
provided that $p\le s\le \frac{2N}{N-4}$ for $N\ge 5$ and
$p\le s <+\infty$ for $N=2,3,4$. Similarly, from (4.3), (4.4), it follows that
\begin{gather*}
S^i_k f_n\text{ and }S^i_k f\in L^s,
\\
\left|S^i_kf_n-S^i_kf\right|_{L^s}\rightarrow 0\text{ as
}n\rightarrow\infty;
\\
S^{\Delta}_kf_n\text{ and }S^{\Delta}_k f\in L^s,
\\
\left|S^{\Delta}_kf_n-S^{\Delta}_kf\right|_{L^s}\rightarrow 0\text{ as }n\rightarrow\infty,
\end{gather*}
provided that $p\le s\le 2$.

By Theorem \ref{Theorem4.1}, we know that $T_k f\in C^4$ and that
 $\partial_i T_kf_n=S^i_kf_n$ for $i=1,2,\dots,N$. Putting $s=2$ in
the preceding statements we deduce that $T_k f\in H^2$ with
$\partial_i T_k f=S^i_kf$ for $i=1,2,\dots,N$,
and $\Delta T_kf=S^{\Delta}_kf$. Furthermore, setting $w=T_kf$ for any
$v\in C^{\infty}_0(\R^N)$, we have
\begin{align*}
\int\Delta w\Delta v-(\lambda w+f)v\,dx
&= \int w\Delta^2 v-\left(\lambda w+f\right)v \,dx\\
&= \lim_{n\rightarrow\infty}\int T_k f_n\Delta  v
-\left(\lambda T_k f_n+f_n\right)v \,dx\\
&= \lim_{n\rightarrow\infty}\int \Delta (T_k f_n)
\Delta v-\left(\lambda T_k f_n +f_n\right)v\,dx\\
&= \lim_{n\rightarrow\infty}\int\left(\Delta^2 (T_k f_n)
+k^2 T_k f_n-f_n \right)v \,dx=0
\end{align*}
by Theorem \ref{Theorem4.1}. This proves that $w$ is a weak solution of (3.1).
\end{proof}

Having established this relationship between weak solutions and
convolutions with fundamental solutions,
 we now have a better understanding of the
 regularity of the weak solutions.

\begin{theorem}
\label{Theorem4.3} Let $f\in L^p\cap L^q$ where $p\in\left(\frac{2N}{N+4},
2\right]$ and $q\ge p$.
Let $u$ be a solution of \textup{(3.1)} for $f$ and some $\lambda\in {\R}$. Then
\begin{enumerate}
\item[{\textup{i)}}]   $u\in W^{2,s} (\R^N)$ \rlap{where}
\[
\begin{array}{ll}
p\le s\le \infty & {\text{ if }}  q>\frac{N}{2}\,,\\
p\le s<\infty & {\text{ if }}  q=\frac{N}{2}\,,\\
p\le s<\frac{Nq}{N-2q} & {\text{ if }}  q<\frac{N}{2}\,;
\end{array}
\]

\item[{\textup{ii)}}]  if $q>\frac{N}{4}$, then $u\in L^{\infty}\cap  C$ and
\[
\lim_{\left|x\right|\rightarrow\infty}u(x)=0;
\]

\item[{\textup{iii)}}]  $u\in L^s$ where
\[
\begin{array}{ll}
p\le s<\infty & {\text{ if }}  q=\frac{N}{4}\,,\\
p\le s<\frac{Nq}{N-4q} & {\text{ if }}  q<\frac{N}{4}\,.
\end{array}
\]
\end{enumerate}
\end{theorem}

\begin{proof}
{\textup{i)}} \ For all $v\in C^{\infty}_0(\R^N)$,
$$
\int\Delta u \Delta v \,dx
=\int\left(\lambda u+f\right)v\,dx
=\int\left(-u+g\right)v\,dx\,,
$$
where $g=(\lambda+1)u+f$. Now since $u\in H^2(\R^N)$, $(\lambda+1)u\in L^r$ for $2\le r<+\infty$ if $N=2,3,4$ and
$2\le r<\frac{2N}{N-4}$ for $N\ge 5$. Thus
$u$ is a weak solution of
$$\Delta ^2 u=-u+g,$$
and so $u=T_1(\lambda+1)u+T_1 f$ (from Theorem \ref{Theorem4.2}). Using Lemma \ref{Lemma3.3} and
a bootstrap argument to deal with the term $T_1(\lambda+1)u$,
the result now follows from (4.1) to (4.6).

{\textup{ii)}} \ By {\textup{i)}}, $u\in W^{2,s}$ for some $s>\frac{N}{2}$ provided
that $q>\frac{N}{4}$. For $s>\frac N2$, $W^{2,s}\hookrightarrow C\cap L^{\infty}$ and
$\lim\limits_{\left|x\right|\rightarrow\infty}u(x)=0$ for  all $u\in W^{2,s}$.

{\textup{iii)}} \ This follows from {\textup{i)}} and the Sobolev inclusions, or
 directly from (4.1).
\end{proof}

\section{\label{S5Reg}Regularity for nonlinear equations}

In this section, we establish the regularity of weak solutions of (1.1).

\begin{proof}[ Proof of Theorem \textup{\ref{Theorem1.1}}]
  Let  $u(x)$ be
a
solution of (1.1). Set
$f=g(x, u(x))-a(x) u(x)$. Then $u(x)$ must be
a
solution of
$$
\left\{
\begin{array}{lr}
\Delta^2 u=f,\\
u\in H^2 (\R^N).
\end{array}
\right.
\eqno{(5.1)}
$$

{}From the assumptions $H_2)$ and $H_3)$ and $1\le \sigma+1<\frac{N+4}{N-4}$,
it follows that for all $u\in H^2(\R^N)$,
$$
f\in L^p \text{\qquad where }
\left\{
\begin{array}{ll}
2\le(\sigma+1)p<\infty &{\text{ if }} N=2,3,4,\\
2\le(\sigma+1)p\le \frac{2N}{N-4} &{\text{ if }} N\ge 5.
\end{array}
\right.
$$
Now
$$
\left[\frac{2}{\sigma+1},\frac{2N}{N-4}\cdot\frac{1}{\sigma+1}\right]
\cap \left(\frac{2N}{N+4},2\right]\not=\emptyset,
$$
since the restrictions on $\sigma$ ensure that
$\frac{2}{\sigma+1}\le 2$ and $\frac{2N}{N-4}\cdot\frac{1}{\sigma+1}>\frac{2N}{N+4}$.
Thus $u$ is a weak solution of (3.1) for $f$ and $\lambda=0$ where $f\in L^p$ for some
$p\in\left(\frac{2N}{N+4},2\right]$. From Theorem \ref{Theorem4.3},
$u\in W^{2,s}(\R^N)$ for
some $s>2$, and so
$$
f\in L^p \text{\qquad where }
\left\{
\begin{array}{ll}
2\le (\sigma+1)p\le\infty &{\text{ if }}s>\frac{N}{2}\,,\\
2\le (\sigma+1)p<\infty &{\text{ if }} s=\frac{N}{2}\,,\\
2\le(\sigma+1)p\le\frac{Ns}{N-2s} &{\text{ if }}  s<\frac{N}{2}\,.
\end{array}
\right.
$$
Noting that $\frac{Ns}{N-2s}>\frac{2N}{N-4}$ for $2<s<\frac{N}{2}$, we see by
Lemma \ref{Lemma3.2} and
a bootstrap argument
that
$u\in W^{2,s}$ for all $2\le s\le \infty$.
This implies that $u\in L^{\infty}$ and so $f\in L^2$. Again by Lemma \ref{Lemma3.2} we now
 also have $u\in H^4$.
\end{proof}

\renewcommand{\footnotesize}{\normalsize}

\end{document}